\documentclass[11pt,twoside]{article}
\usepackage{amsfonts}
\usepackage[plainpages=false]{hyperref}
\usepackage{amsfonts,latexsym,rawfonts,amsmath,amssymb,amsthm}
\usepackage{amsmath,amssymb,amsfonts,latexsym,lscape,rawfonts,appendix}
\newcommand{\pd}[2]{\frac {\partial #1}{\partial #2}}

\newcommand{\oo}{\omega}

\newcommand{\ee}{\epsilon}

\newcommand{\beq}{\begin{equation}}
\newcommand{\eeq}{\end{equation}}
\newcommand{\beqs}{\begin{eqnarray*}}
\newcommand{\eeqs}{\end{eqnarray*}}
\newcommand{\beqn}{\begin{eqnarray}}
\newcommand{\eeqn}{\end{eqnarray}}
\newcommand{\beqa}{\begin{array}}
\newcommand{\eeqa}{\end{array}}

\def\lra{\longrightarrow}

\def\bc{\begin{center}}
\def\ec{\end{center}}

\def\cP{{\cal P}}

\def\cF{{\mathcal F}}

\def\cP{{\mathcal P}}

\def\RR{{\mathbb R}}

\def\begeq{\begin{equation}}
\def\endeq{\end{equation}}
\def\and{\quad{\rm and}\quad}

\let\lra=\longrightarrow

\def\mapright\#1{\,\smash{\mathop{\lra}\limits^{\#1}}\,}

\def\pbp{\sqrt{-1}\partial\bar\partial}

\def\an{\;\;\;{\rm and}\;\;\;}
\newtheorem{prop}{Proposition}[section]
\newtheorem{theo}[prop]{Theorem}
\newtheorem{lem}[prop]{Lemma}
\newtheorem{claim}[prop]{Claim}
\newtheorem{cor}[prop]{Corollary}
\newtheorem{rem}[prop]{Remark}

\newtheorem{defi}[prop]{Definition}

\title{{\bf\Large{A new formula for the energy functionals $E_k$ and its applications}}}

\author{Haozhao Li }

\begin{document}
\bibliographystyle{plain}
\maketitle
\begin{abstract}We give a new formula for the energy functionals $E_k$
defined by Chen-Tian \cite{[chen-tian1]}, and discuss the relations
between these functionals. We also apply our formula to give a new
proof of the fact that the holomorphic invariants corresponding to
the $E_k$ functionals are equal to the Futaki invariant.
\end{abstract}
\section{Introduction}
In \cite{[chen-tian1]}, a series of energy functionals $E_k(k=0, 1,
\cdots, n)$ were introduced by X.X. Chen and  G. Tian which were
used to prove the convergence of the K\"ahler Ricci flow under some
curvature assumptions. The first energy functional $E_0$ of this
series is exactly the $K$-energy introduced by Mabuchi in
\cite{[Ma]}, which can be defined for any K\"ahler potential
$\varphi(t)$ on a K\"ahler manifold $(M, \oo)$ as follows:
$$\frac {d}{dt}E_0(\varphi(t))=-\frac 1V\int_M\;\pd {\varphi}t(R_{\varphi}-r)\oo_{\varphi}^n.$$
Here $R_{\varphi}$ is the scalar curvature with respect to the
K\"ahler metric $\oo_{\varphi}=\oo+\pbp \varphi$,   $r=\frac
{[c_1(M)]\dot[\oo]^{n-1}}{[\oo]^n}$ is the average of $R_{\varphi}$
and $V=[\oo]^n$ is the volume.

It is well-known that the behavior of the $K$-energy plays a central
role on the existence of K\"ahler-Einstein metrics and constant
scalar curvature metrics. In \cite{[BaMa]}, Bando-Mabuchi proved
that the $K$-energy is bounded from below on a K\"ahler-Einstein
manifold with $c_1(M)>0$. It has been shown by G. Tian in
\cite{[Tian2]}\cite{[Tian3]} that $M$ admits a K\"ahler-Einstein
metric if and only if  the $K$-energy is proper. Recently, Chen-Tian
in \cite{[chen-tian3]} extended these results to extremal K\"ahler
metrics, and Cao-Tian-Zhu in \cite{[CTZ]}\cite{[Tian-Zhu]} proved
similar results on  K\"ahler Ricci solitons. So a natural question
is how the energy functionals $E_k $ are related to these extremal
metrics.

Following a question posed by Chen in  \cite{[chen1]}, Song-Weinkove
recently  proved  in \cite{[SoWe]}  that the energy functionals
$E_k$ have a lower bound  on the space of K\"ahler metrics with
nonnegative Ricci curvature for K\"ahler-Einstein manifolds.
Moreover, they also showed that modulo holomorphic vector fields,
$E_1$
  is proper if and only if there exists a K\"ahler-Einstein metric. Shortly afterwards, N. Pali
    \cite{[Pali]} gave a formula between $E_1$ and the $K$-energy $E_0$, which
    implies $E_1$ has a lower bound if the $K$-energy is bounded from below.
    Tosatti  \cite{[Tosatti]} proved under some curvature assumptions,
    the critical point of $E_k$ is a K\"ahler-Einstein metric.   Pali's theorem says that the functional $E_1$
    is always bigger than the $K$-energy. However, we proved that the converse
    is also true in \cite{[CLW]}. Following suggestion of X. X. Chen, we set out to investigate
    the relations  between these energy
    functionals for the general case; in particular, the relations about lower bounds of these functionals.\\

Now we state our results. Let $M $ be an $n$-dimensional compact
K\"ahler manifold with  $c_1(M)>0$, and $\oo$ be a fixed K\"ahler
metric in the K\"ahler class $2\pi c_1(M).$ Write
$$\cP(M,  \oo)=\{\varphi\in C^{\infty}(M, \RR)\;|\; \oo_{\varphi}=\oo+\pbp \varphi>0 \;\;{\rm on}\;\; M\}.$$
For any $k=0, 1,\cdots, n $, we define the functional $E_{k,
\oo}^0(\varphi)$  on $\cP(M, \oo)$ by
$$E_{k, \oo}^0(\varphi)=\frac 1V\int_M\; \Big(\log \frac {\oo^n_{\varphi}}{\oo^n}-h_{\oo}\Big)\Big(
\sum_{i=0}^k\; Ric_{\varphi}^i\wedge \oo^{k-i}\Big)\wedge
\oo_{\varphi}^{n-k}+\frac 1V\int_M\; h_{\oo} \Big(\sum_{i=0}^k\;
Ric_{\oo}^i\wedge \oo^{k-i}\Big)\wedge \oo^{n-k}.$$Here $h_{\oo}$
is the Ricci potential defined by
$$Ric_{\oo}-\oo=\pbp h_{\oo}, \an \int_M\;(e^{h_{\oo}}-1)\oo^n=0.$$
Let $\varphi(t)(t\in [0, 1])$ be a path from $0$ to $\varphi$ in
$\cP(M, \oo)$, we define
$$J_{k, \oo}(\varphi)=-\frac {n-k}{V}\int_0^1\;\int_M\;
\pd {\varphi(t)}{t}(\oo_{\varphi(t)}^{k+1}-\oo^{k+1})\wedge
\oo_{\varphi(t)}^{n-k-1}\wedge dt.$$ Then the functional $E_{k,
\oo}$ is defined as follows
$$E_{k, \oo}(\varphi)=E_{k, \oo}^0(\varphi)-J_{k, \oo}(\varphi).$$
For simplicity, we will often drop the subscript $\oo$ and write
$E_k$ instead of $E_{k, \oo}(\varphi).$ The main result of this
paper is the following
\begin{theo}\label{theo1}For any $k=1, 2, \cdots, n $, we have
$$\sum_{i=0}^k\;(-1)^i\binom{k+1}{i+1}E_{i, \oo}(\varphi)=\frac 1V\int_M\; u(\pbp u)^k\wedge \oo_{{\varphi}}^{n-k}
+\frac 1V\int_M\; h_{\oo}(-\pbp h_{\oo})^k\wedge \oo^{n-k},$$
where
$$u=\log \frac {\oo^n_{\varphi}}{\oo^n}+\varphi-h_{\oo}.$$
\end{theo}

\begin{rem}Theorem \ref{theo1} generalizes Pali's formula in
\cite{[Pali]}. In fact, when $k=1, 2$, we have the following \beqs
2E_0-E_1&=&-\frac 1V \int_M\;\sqrt{-1}\partial u\wedge \bar\partial
u\wedge \oo_{\varphi}^{n-1}+c_1,\\3E_0-3E_1+E_2&=&-\frac 1V
\int_M\;\sqrt{-1}\partial u\wedge \bar\partial u\wedge\pbp u\wedge
\oo_{\varphi}^{n-2}+c_2, \eeqs where $c_1, c_2$ are two constants
depending only on $\oo$.

\end{rem}

Next we use Theorem \ref{theo1} to get the lower bound of $E_k.$
\begin{theo}\label{theo2}For any positive integer $k=2, \cdots, n$, and any K\"ahler metric
  $\oo_{\varphi}$ satisfying $Ric_{\varphi}\geq -\frac 2{k-1}\oo_{\varphi}$, we have
$$E_k(\varphi)\geq (k+1)E_0(\varphi)+c_k,$$
where $c_k$ is a constant defined by \beq c_k=\frac
1V\int_M\;\sum_{i=0}^{k-1}\;(-1)^{k-i}\binom{k+1}{i}h_{\oo} (-\pbp
h_{\oo})^{k-i}\wedge \oo^{n-k+i}. \label{ck}\eeq
\end{theo}

\begin{rem}Theorem \ref{theo2} generalizes some of Song-Weinkove's results in
\cite{[SoWe]}. Since $E_0$ is bounded from below on $\cP(M, \oo)$ on
a K\"ahler-Einstein manifold, from Theorem \ref{theo2} we obtain
lower bounds on the functionals $E_k$ under some weaker conditions.
\end{rem}

\begin{rem}In \cite{[CLW]}, we proved that $E_1$ is bounded from below
if and only if $E_0$ is bounded from below on $\cP(M, \oo)$. Using
the same method, we also prove that $E_0$ is bounded from below if
and only if the $F$ functional defined by Ding-Tian \cite{[DiTi]} is
bounded from below in \cite{[Li]}. We expect that the lower
boundedness of these functionals are equivalent on $\cP(M, \oo)$ in
\cite{[CLW]}.
\end{rem}

Finally, we will prove that all the Chen-Tian holomorphic invariants
$\cF_k$ defined by $E_k$ are the Futaki invariant in the canonical
K\"ahler class.

\begin{theo}\label{theo3}For all $k=0,1,\cdots, n$, we have
$$\cF_k(X, \oo)=(k+1)\cF_0(X, \oo).$$
\end{theo}

\begin{rem}This result was first proved by C. Liu in \cite{[Liu]},
and here  we give a new proof by using our formula. However, these
two methods are essentially the same.
\end{rem}

\noindent {\bf Acknowledgements}:  This work was done while I was
attending the summmer school on geometric analysis in University of
Science and Technology of China (USTC) in 2006, and I would like to
express thanks to USTC. I would also like to thank Professor X. X.
Chen, W. Y. Ding and X. H. Zhu for their constant support and
advice. Thanks also go to Y. Rubinstein, V. Tosatti for pointing out
some mistakes in Theorem \ref{theo2},  B. Wang , W. Y. He for
carefully reading the draft, and the referees for  numerous
suggestions which helped to improve the presentation.

\section{A new formula on $E_k$}
In this section, we will prove Theorem \ref{theo1} and Corollary
\ref{maincor1}.\\
{\it Proof of Theorem \ref{theo1}.} By the definition of $u$, we
have
$$\pbp u=-Ric_{\varphi}+\oo_{\varphi}.$$
Therefore, we have \beqs \Big( \sum_{p=0}^i\; Ric_{\varphi}^p\wedge
\oo^{i-p}\Big)\wedge \oo_{\varphi}^{n-i}=\Big( \sum_{p=0}^i\;
(\oo_{\varphi}-\pbp u)^p\wedge (\oo_{\varphi}-\pbp
\varphi)^{i-p}\Big)\wedge \oo_{\varphi}^{n-i} \eeqs By the
definition of $E_k^0$ we have \beqs &
 &\sum_{i=0}^k\;(-1)^i\binom{k+1}{i+1}E_{i}^0(\varphi)\\&=&\frac 1V\int_M\;
(u-\varphi)\Big(
\sum_{i=0}^k\;(-1)^i\binom{k+1}{i+1}\sum_{p=0}^i\;
(\oo_{\varphi}-\pbp u)^p\wedge (\oo_{\varphi}-\pbp
\varphi)^{i-p}\Big)\wedge \oo_{\varphi}^{n-i}\\
&&+\frac 1V\int_M\; h_{\oo}\Big(\sum_{i=0}^k\;(-1)^i\binom
{k+1}{i+1}\sum_{p=0}^i\;(\oo+\pbp h_{\oo})^p\wedge
\oo^{i-p}\Big)\wedge \oo^{n-i}. \eeqs Now we have the following
lemma:
\begin{lem}For any two variables $x,  y$ and any integer $k>0$, we
have
\begin{enumerate}
    \item
    \beq\sum_{i=0}^k(-1)^i\binom{k+1}{i+1}\sum_{p=0}^{i}\;(1-x)^p(1-y)^{i-p}=\sum_{i=0}^k\;x^{k-i}y^{i},\label{a1}\eeq
    \item
    \beq\sum_{i=0}^k(-1)^i\binom{k+1}{i+1}\sum_{p=0}^{i}\;(1+x)^p=(-x)^k.\label{a2}\eeq
\end{enumerate}

\end{lem}
\begin{proof} By direct calculation, we have
\beqs &
&(x-y)\sum_{p=0}^k(-1)^p\binom{k+1}{p+1}\sum_{i=0}^{p}\;(1-x)^i(1-y)^{p-i}\\&=&
\sum_{p=0}^k
\binom{k+1}{p+1}((x-1)^{p+1}-(y-1)^{p+1})\\&=&x^{k+1}-y^{k+1}. \eeqs
Then the equality (\ref{a1}) holds. Similarly, we can prove the
equality (\ref{a2}).
\end{proof}

Thus, the energy functionals $E_k^0$ satisfy the equality
\beqn\sum_{i=0}^k\;(-1)^i\binom{k+1}{i+1}E_{i}^0(\varphi)&=&\sum_{i=0}^k\;\frac
1V\int_M\;(u-\varphi)(\pbp u)^{k-i}\wedge (\pbp \varphi )^{i}\wedge
\oo_{\varphi}^{n-k}\nonumber\\&& +\frac 1V\int_M\; h_{\oo}(-\pbp
h_{\oo})^k\wedge \oo^{n-k}.\label{eq3}\eeqn Observe that for $0\leq
i\leq k-1,$ \beqs & & \int_M\;(u-\varphi)(\pbp u)^{k-i}\wedge (\pbp
\varphi )^{i}\wedge
\oo_{\varphi}^{n-k}\\&=&\int_M\;u\pbp(u-\varphi)\wedge(\pbp
u)^{k-i-1}\wedge (\pbp \varphi )^{i}\wedge \oo_{\varphi}^{n-k}\\
&=&\int_M\;u(\pbp u)^{k-i}\wedge (\pbp \varphi )^{i}\wedge
\oo_{\varphi}^{n-k}\\
& &-\int_M\;u(\pbp u)^{k-i-1}\wedge (\pbp \varphi )^{i+1}\wedge
\oo_{\varphi}^{n-k}. \eeqs Thus, the equality (\ref{eq3}) can be
written as \beqn
\sum_{i=0}^k\;(-1)^i\binom{k+1}{i+1}E_{i}^0(\varphi)&=&\frac
1V\int_M\;u(\pbp u)^{k}\wedge \oo_{\varphi}^{n-k}-\frac
1V\int_M\;\varphi(\pbp \varphi)^{k}\wedge
\oo_{\varphi}^{n-k}\nonumber\\&&+\frac 1V\int_M\; h_{\oo}(-\pbp
h_{\oo})^k\wedge \oo^{n-k}.\label{eq4}\eeqn

Next we calculate $J_k(\varphi)$ via a linear path $t\varphi\in
\cP(M, \oo)$ for $t\in [0, 1].$ By the definition of $J_k$ we have
\beqs & &\sum_{i=0}^k\;(-1)^i\binom{k+1}{i+1}J_i(\varphi)\\
&=&\frac
1V\int_0^1\;\int_M\;\sum_{i=0}^k\;-(n-i)(-1)^i\binom{k+1}{i+1}\varphi(\oo_{t\varphi}^{i+1}-(\oo_{t\varphi}-t\pbp\varphi)^{i+1})\wedge
\oo_{t\varphi}^{n-i-1}\wedge dt. \eeqs It is easy to check the
following lemma:
\begin{lem}
Let $B_i=-(n-i)(1-(1-x)^{i+1})$, for any integer $k\geq 1$ we have
$$\sum_{i=0}^k\;(-1)^i\binom {k+1}{i+1}B_i=-(n-k)x^{k+1}-(k+1)x^k.$$
\end{lem}

Thus, we have \beqs & &\sum_{i=0}^k\;(-1)^i\binom{k+1}{i+1}J_i(\varphi)\\
&=& \frac 1V\int_0^1\;\int_M\;-(n-k)t^{k+1}\varphi
(\pbp\varphi)^{k+1}\wedge \oo_{t\varphi}^{n-k-1}\wedge dt\\& &-
\frac 1V\int_0^1\;\int_M\;(k+1)t^{k}\varphi (\pbp\varphi)^{k}\wedge
\oo_{t\varphi}^{n-k}\wedge dt\\&=& \frac 1V\int_0^1\;\int_M\;-\frac
{d}{dt}\Big(t^{k+1}\varphi (\pbp\varphi)^{k}\wedge
\oo_{t\varphi}^{n-k}\Big)\wedge dt\\
&=&-\frac 1V\int_M\;\varphi(\pbp \varphi)^{k}\wedge
\oo_{\varphi}^{n-k}.\eeqs Combining this with the equality
(\ref{eq4}), we have
$$\sum_{i=0}^k\;(-1)^i\binom{k+1}{i+1}E_{i}(\varphi)=\frac 1V\int_M\; u(\pbp u)^k\wedge \oo_{{\varphi}}^{n-k}+\frac 1V\int_M\;
h_{\oo}(-\pbp h_{\oo})^k\wedge \oo^{n-k}.$$ \hfill $\square$

Next we will use Theorem \ref{theo1} to prove the following
corollary.
\begin{cor}\label{maincor1}Let $$F_k(\varphi)=\frac 1V\int_M\; u(\pbp u)^k\wedge \oo_{{\varphi}}^{n-k}+\frac 1V\int_M\;
h_{\oo}(-\pbp h_{\oo})^k\wedge \oo^{n-k},$$ we have
\begin{enumerate}

  \item For nonnegative integers $p, k\;(0\leq p\leq k-2\leq n-2)$, we have
\beq\sum_{i=p}^k\;(-1)^i\binom{k-p}{i-p}E_i=\sum_{i=0}^{p+1}\;(-1)^i\binom{p+1}{i}F_{k-i}.\label{c1}\eeq

\item For any positive integer $k=1, 2, \cdots, n$, we have \beq
E_k-E_{k-1}-E_0=\frac 1V\int_M\;
u\Big(Ric_{\varphi}^k-\oo_{\varphi}^k\Big)\wedge
\oo_{\varphi}^{n-k}+\frac 1V\int_M\;
h_{\oo}\Big(Ric_{\oo}^k-\oo^k\Big)\wedge \oo^{n-k}.\label{c2}\eeq
  \item For any positive integer $k=1, 2, \cdots, n$, we have
\beq
E_k=\sum_{i=0}^{k-1}\;(-1)^{k-i}\binom{k+1}{i}F_{k-i}+(k+1)E_0.\label{c3}\eeq
\end{enumerate}
\end{cor}

\begin{proof} (1). We show this by induction on $p.$ The corollary holds for $p=0.$ In fact, by Theorem \ref{theo1} we have \beqn
\sum_{i=0}^k\;(-1)^i\binom{k+1}{i+1}E_{i}&=&F_k,\label{eq5}\\\sum_{i=0}^{k-1}\;(-1)^i\binom{k}{i+1}E_{i}&=&F_{k-1}.\label{eq6}\eeqn
 Subtract (\ref{eq6}) from (\ref{eq5}), we have
$$\sum_{i=p}^k\;(-1)^i\binom{k}{i}E_i=F_k-F_{k-1}.$$
We assume that  the corollary holds for  $p,$ then \beqn
\sum_{i=p}^k\;(-1)^i\binom{k-p}{i-p}E_{i}&=&\sum_{i=0}^{p+1}\;(-1)^i\binom{p+1}{i}F_{k-i}
,\label{eq7}\\\sum_{i=p}^{k-1}\;(-1)^i\binom{k-p-1}{i-p}E_{i}&=&\sum_{i=0}^{p+1}\;(-1)^i\binom{p+1}{i}F_{k-i-1}
=\sum_{i=1}^{p+2}\;(-1)^{i-1}\binom{p+1}{i-1}F_{k-i}
.\label{eq8}\eeqn  Subtract (\ref{eq8}) from (\ref{eq7}), we have
\beqs
\sum_{i=p+1}^{k}\;(-1)^i\binom{k-p-1}{i-p-1}E_{i}&=&\sum_{i=0}^{p+1}\;(-1)^i\binom{p+1}{i}F_{k-i}
-\sum_{i=1}^{p+2}\;(-1)^{i-1}\binom{p+1}{i-1}F_{k-i}\\
&=&F_k+\sum_{i=1}^{p+1}\;(-1)^i\Big(\binom{p+1}{i}+\binom{p+1}{i-1}\Big)F_{k-i}+(-1)^{p+2}F_{k-p-2}\\
&=&\sum_{i=0}^{p+2}\;(-1)^i\binom{p+2}{i}F_{k-i}. \eeqs The
corollary holds for $p+1$. Thus, the equality (\ref{c1}) holds.

 (2) We can show
the following formula by induction: \beq
E_k-E_{k-1}=\sum_{i=0}^{k-1}\;(-1)^{k-i}\binom{k}{i}F_{k-i}+E_0.\label{t1}\eeq
In fact, by Theorem \ref{theo1} the formula  (\ref{t1}) holds for
$k=1.$  We assume the formula (\ref{t1}) holds for some integer
$k\leq n-1,$ then by (1) we have \beq
E_{k+1}=2E_{k}-E_{k-1}+\sum_{i=0}^{k}\;(-1)^{i-k-1}\binom{k}{i}F_{k+1-i}.\label{t2}\eeq
Thus, we have \beqs
E_{k+1}-E_{k}&=&E_k-E_{k-1}+\sum_{i=0}^{k}\;(-1)^{i-k-1}\binom{k}{i}F_{k+1-i}\\
&=&\sum_{i=0}^{k-1}\;(-1)^{k-i}\binom{k}{i}F_{k-i}+E_0+\sum_{i=0}^{k}\;(-1)^{i-k-1}\binom{k}{i}F_{k+1-i}\\
&=&E_0+\sum_{i=0}^{k}\;(-1)^{k+1-i}\binom{k+1}{i}F_{k+1-i}. \eeqs
Then the formula (\ref{t1}) holds for $k+1$.

On the other hand, by direct calculation we have \beqs
\sum_{i=0}^{k-1}\;(-1)^{k-i}\binom{k}{i}F_{k-i}&=&\frac 1V\int_M\;
u\Big(\sum_{i=0}^{k}\;(-1)^{k-i}\binom{k}{i}(\pbp u)^{k-i}\wedge
\oo_{\varphi}^i-\oo_{\varphi}^k\Big)\wedge \oo_{\varphi}^{n-k}
\\&&+\frac 1V\int_M\;
h_{\oo}\Big(\sum_{i=0}^{k}\;(-1)^{k-i}\binom{k}{i}(-\pbp
h_{\oo})^{k-i}\wedge \oo^i-\oo^k\Big)\wedge \oo^{n-k}\\
&=&\frac 1V\int_M\; u\Big((\oo_{\varphi}-\pbp
u)^k-\oo_{\varphi}^k\Big)\wedge \oo_{\varphi}^{n-k}+\frac
1V\int_M\; h_{\oo}\Big(Ric_{\oo}^k-\oo^k\Big)\wedge \oo^{n-k}\\
&=&\frac 1V\int_M\; u\Big(Ric_{\varphi}^k-\oo_{\varphi}^k\Big)\wedge
\oo_{\varphi}^{n-k}+\frac 1V\int_M\;
h_{\oo}\Big(Ric_{\oo}^k-\oo^k\Big)\wedge \oo^{n-k}. \eeqs Then the
equality (\ref{c2}) holds.

(3). We prove this result by induction on $k.$ The corollary holds
for $k=1$ obviously. We assume that it holds for integers less
than $k,$ then by (1) we have
$$E_k=2E_{k-1}-E_{k-2}+\sum_{i=0}^{k-1}\;(-1)^{i-k}\binom{k-1}{i}F_{k-i}.$$
By induction, we have \beqs
E_{k-1}=\sum_{i=0}^{k-2}\;(-1)^{k-i-1}\binom{k}{i}F_{k-i-1}+kE_0=
\sum_{i=1}^{k-1}\;(-1)^{k-i}\binom{k}{i-1}F_{k-i}+kE_0,\eeqs and
 \beqs
E_{k-2}=\sum_{i=0}^{k-3}\;(-1)^{k-i-2}\binom{k-1}{i}F_{k-i-2}+(k-1)E_0=
\sum_{i=2}^{k-1}\;(-1)^{k-i}\binom{k-1}{i-2}F_{k-i}+(k-1)E_0.\eeqs
Then we have \beqs
E_k&=&2\Big(\sum_{i=1}^{k-1}\;(-1)^{k-i}\binom{k}{i-1}F_{k-i}+kE_0\Big)
-\Big(\sum_{i=2}^{k-1}\;(-1)^{k-i}\binom{k-1}{i-2}F_{k-i}+(k-1)E_0\Big)\\
& &+\sum_{i=0}^{k-1}\;(-1)^{i-k}\binom{k-1}{i}F_{k-i}\\
&=&\sum_{i=0}^{k-1}\;(-1)^{k-i}\binom{k+1}{i}F_{k-i}+(k+1)E_0. \eeqs
Then the equality (\ref{c3}) holds.

\end{proof}

\section{Applications of the new formula}
In this section, we will prove Theorem \ref{theo2} and
\ref{theo3}.
\subsection{On the lower bound of $E_k$}

{\it Proof of Theorem \ref{theo2}.} By the equality (\ref{c3}) of
Corollary \ref{maincor1}, we have \beqs
E_k-(k+1)E_0&=&\sum_{i=0}^{k-1}\;(-1)^{k-i}\binom{k+1}{i}F_{k-i}\\
&=&\frac 1V\int_M\;\sum_{i=0}^{k-1}\;(-1)^{k-i}\binom{k+1}{i}u (\pbp
u)^{k-i}\wedge \oo_{\varphi}^{n-k+i}+c_k\\&=& \frac
1V\int_M\;\sqrt{-1}\partial u\wedge\bar \partial u\wedge
\Big(\sum_{i=0}^{k-1}(-1)^{k-i-1}\binom{k+1}{i}(\pbp
u)^{k-i-1}\wedge \oo_{\varphi}^i\Big)\wedge \oo_{\varphi}^{n-k}+c_k
\\&=& \frac
1V\int_M\;\sqrt{-1}\partial u\wedge\bar \partial u\wedge
\Big(\sum_{i=0}^{k-1}\binom{k+1}{i}(Ric_{\varphi}-\oo_{\varphi})^{k-i-1}\wedge
\oo_{\varphi}^i\Big)\wedge \oo_{\varphi}^{n-k}+c_k, \eeqs where
$c_k$ is a constant defined by (\ref{ck}). Observe that \beq
\sum_{i=0}^{k-1}\binom{k+1}{i}(Ric_{\varphi}-\oo_{\varphi})^{k-i-1}\wedge
\oo_{\varphi}^i=\sum_{i=1}^{k}i Ric_{\varphi}^{k-i}\wedge
\oo_{\varphi}^{i-1}.\label{x1}\eeq Then we need to check when
(\ref{x1}) is nonnegative. Obviously, this is true when
$Ric_{\varphi}\geq 0,$ Here we want to get a better condition on
Ricci curvature. If $k=2$, we need to assume $Ric_{\varphi}\geq -2
\oo_{\varphi}.$ Now we assume $k\geq 3.$ Set
$$P(x)=\sum_{i=1}^{k}ix^{k-i}=(x+\frac
2{k-1})^{k-1}+\sum_{i=2}^k\;a_i (x+\frac 2{k-1})^{k-i},$$ where
$a_i$ are the constants defined by
$$a_i=\frac 1{(k-i)!}P^{(k-i)}(-\frac 2{k-1}).$$
By Lemma {A.1} in the appendix, $a_i\geq 0.$ Then if
$Ric_{\varphi}\geq -\frac 2{k-1}\oo_{\varphi},$ we have
$$\sum_{i=1}^{k}i Ric_{\varphi}^{k-i}\wedge
\oo_{\varphi}^{i-1}=\Big(Ric_{\varphi}+\frac
2{k-1}\oo_{\varphi}\Big)^{k-1}+\sum_{i=2}^k\;a_i
\Big(Ric_{\varphi}+\frac 2{k-1}\oo_{\varphi}\Big)^{k-i}\wedge
\oo_{\varphi}^{i-1}\geq 0.$$ Therefore, $E_k\geq (k+1)E_0+c_k.$
\hfill $\square$

\subsection{On the holomorphic invariants $\cF_k$}
In this subsection, we will use the equality (\ref{c3}) of Corollary
\ref{maincor1} to prove that all the holomorphic invariants defined
in \cite{[chen-tian1]} are the Futaki invariant. This result was
first obtained by Liu in \cite{[Liu]}. Here we give a new proof by
using our formula.

Let $X$ be a holomorphic vector field. Then by $c_1(M)>0$, we can
decompose $i_X\oo$ as $i_X\oo=\sqrt{-1}\bar\partial \theta_X,$ where
$\theta_X$ is a potential function of $X$ with respect to $\oo$.
\begin{defi}(cf. \cite{[chen-tian1]}) For any holomorphic vector field $X$, we define
$$\cF_k=(n-k)\int_M\;\theta_X \oo^n+\int_M\; \Big((k+1)\Delta \theta_X Ric_{\oo}^k\wedge \oo^{n-k}-(n-k)\theta_X
Ric_{\oo}^{k+1}\wedge \oo^{n-k-1}\Big).$$
\end{defi}
It was proved that $\cF_k$ is a holomorphic invariant. When $k=0,$
we have
$$\cF_0(X, \oo)=n\int_M\; X(h_{\oo}) \oo^n,$$
which is  a multiple of the Futaki invariant.
\begin{prop}\label{prop}(cf. \cite{[chen-tian1]}) Let $\{\Phi(t)\}_{|t|<\infty}$ be the one-parameter
subgroup of automorphisms induced by $Re(X)$. Then
$$\frac {dE_k(\varphi_t)}{dt}=\frac 1VRe(\cF_k(X, \oo)),$$
where $\varphi_t$ are the K\"ahler potentials of $\Phi_t^*\oo$,
i.e.,$\Phi_t^*\oo=\oo+\pbp \varphi_t.$
\end{prop}

Now we can prove Theorem \ref{theo3}. \\
{\it Proof of Theorem
\ref{theo3}.} By Corollary \ref{maincor1},  we only need to show
$$\frac {dF_k(\varphi_t)}{dt}=0,$$
for all $k$, where $\varphi_t$ is the K\"ahler potential defined in
the previous proposition. Differentiating
$\oo_{\varphi}=\Phi_t^*\oo=\oo+\pbp \varphi_t,$ we get
$$L_{Re(X)}\oo_{\varphi}=\pbp \pd {\varphi_t}{t}.$$
On the other hand, since $L_X\oo=\pbp \theta_X,$ we have
$$\pd {\varphi_t}t=Re(\theta_X(\varphi))+c,$$
where $c$ is a constant and
$\theta_X(\varphi)=\theta_X+X(\varphi).$ By the definition of $u$,
we have
$$Ric_{\varphi}-\oo_{\varphi}=-\pbp u.$$
Take the inner product on both sides, we have
$$-\Delta \theta_X(\varphi)-\theta_X(\varphi)=-X(u).$$
Here $\Delta$ is the Lapacian with respect to $\oo_{\varphi}.$ On
the other hand
$$\pd ut=\Delta \pd {\varphi}t+\pd {\varphi}t=Re(\Delta \theta_X(\varphi)+\theta_X(\varphi))+c=Re(X(u))+c.$$
Thus, \beqs & &\pd {}t\int_M\; u(\pbp u)^k\wedge
\oo_{\varphi}^{n-k}\\&=&\int_M\;\pd ut(\pbp u)^k\wedge
\oo_{\varphi}^{n-k}+\int_M\;ku \pbp \pd ut\wedge (\pbp
u)^{k-1}\wedge \oo_{\varphi}^{n-k}\\& &+\int_M\;(n-k)u(\pbp
u)^k\wedge \pbp \pd {\varphi}t\wedge \oo_{\varphi}^{n-k-1}\\
&=&Re\Big(\int_M\; (k+1)X(u)(\pbp u)^k\wedge
\oo_{\varphi}^{n-k}+(n-k)\theta_X(\varphi) (\pbp u)^{k+1}\wedge
\oo_{\varphi}^{n-k-1}\Big)\\
&=&Re\Big(\int_M\; i_X(\partial u (\pbp u)^k\wedge
\oo_{\varphi}^{n-k})\Big)\\
&=&0. \eeqs Thus, by the equality (\ref{c3}) in Corollary
\ref{maincor1},  we have
$$\frac {dE_k(\varphi_t)}{dt}=\sum_{i=0}^{k-1}\;(-1)^{k-i}
\binom{k+1}{i}\frac {d}{dt}F_{k-i}(\varphi_t)+(k+1)\frac
{dE_0(\varphi_t)}{dt}=\frac {k+1}VRe(\cF_0(X, \oo)).$$ By
Proposition \ref{prop}, the theorem is proved. \hfill $\square$

\appendix


\section{An elementary lemma} In the proof of Theorem \ref{theo2}, we
need to use the following lemma.

\begin{lem}Let  $m$ be a positive integer. Consider the polynomial
$$P(x)=x^m+2x^{m-1}+\cdots+ mx+(m+1),$$
then for any $i(0\leq i\leq m),$ the $i^{th}$ derivative of the
polynomial at the point $x=-\frac 2m$ is nonnegative.
\end{lem}
\begin{proof}The $i^{th}$ derivative of the
polynomial is
$$P^{(i)}(x)=\sum_{p=0}^{m-i}\; (m-i+1-p)(i+p)(i+p-1)\cdots(p+1)x^p.$$
For simplicity, we define $a(p, i)$ by
$$a(p, i)=(i+p)(i+p-1)\cdots(p+1).$$
If $m-i$ is even, then \beq P^{(i)}(x)=a(m-i,
i)x^{m-i}+\sum_{p=0}^{\frac {m-i}2-1} \Big((m-i+1-2p)\; a(2p,
i)x^{2p}+(m-i-2p)a(2p+1,  i)x^{2p+1}\Big).\label{aq0}\eeq If $m-i$
is odd, we write $P^{(i)}(x)$ as \beq P^{(i)}(x)=\sum_{p=0}^{\frac
{m-i-1}2} \Big((m-i+1-2p)\; a(2p,  i)x^{2p}+(m-i-2p)a(2p+1,
i)x^{2p+1}\Big).\label{aq1}\eeq Note that $P^{(m-1)}(-\frac 2m)=0,$
so we can assume $i\leq m-2.$ Since the lemma is trivial for $1\leq
m\leq 10$, we assume $m> 10.$ For simplicity, we define
$$A_p(x)=(m-i+1-2p)\; a(2p,  i)x^{2p}+(m-i-2p)a(2p+1,
i)x^{2p+1}.$$

\begin{claim}\label{claim1}If $1\leq p\leq \frac {m-i-1}{2},$ we have $A_p(-\frac 2m)> 0.$
\end{claim}
\begin{proof}We need to show
$$\frac {(m-i+1-2p)}{m-i-2p}\frac {m(2p+1)}{2(i+2p+1)}>1.$$
Since $1\leq p\leq \frac {m-i-1}{2}, $ this is obvious because
$$\frac {m(2p+1)}{2(i+2p+1)}\geq \frac {3m}{2m}>1.$$
The claim is proved.
\end{proof}

By Claim \ref{claim1}, all the terms on the right hand side of
(\ref{aq0}) and (\ref{aq1}) are positive except  $A_0(-\frac 2m).$
Note that if $0\leq i\leq\frac m2,$
\beqs A_0(-\frac 2m)&=&(m-i)(i+1)\cdots 2(-\frac 2m)+(m-i+1)i(i-1)\cdots 1\\
&=&\frac {i!}m (m-2i)(m-i-1)\\&\geq & 0. \eeqs So it only remains to
deal with the case $i>\frac m2.$ Now, we consider the case $\frac
12m< i\leq m-5.$ The following claim shows that $A_0+A_1+A_2$ is
positive at $x=-\frac 2m$ in this case.

\begin{claim}\label{claim3}If $\frac 12m< i\leq m-5$, then $(A_0+A_1+A_2)(-\frac 2m)>0.$
\end{claim}
\begin{proof}In fact,
\beqs &&\frac {120}{i!}(A_0+A_1+A_2)(-\frac 2m)\\&=&-\frac {32}{m^5}(m-i-4)(i+5)(i+4)(i+3)(i+2)(i+1)+\frac {80}{m^4}(m-i-3)(i+4)(i+3)(i+2)(i+1)\\
&&-\frac {160}{m^3}(m-i-2)(i+3)(i+2)(i+1)+\frac
{240}{m^2}(m-i-1)(i+2)(i+1)\\&&-\frac {240}{m}(m-i)(i+1)+120(m-i+1).
\eeqs Observe that
$$\frac {32}{m^4}(m-i-3)(i+4)(i+3)(i+2)(i+1)>\frac {32}{m^5}(m-i-4)(i+5)(i+4)(i+3)(i+2)(i+1),$$
so we only need to show \beqs A&:=&\frac
{48}{m^4}(m-i-3)(i+4)(i+3)(i+2)(i+1)-\frac
{160}{m^3}(m-i-2)(i+3)(i+2)(i+1)\\&&+\frac
{240}{m^2}(m-i-1)(i+2)(i+1)-\frac {240}{m}(m-i)(i+1)+120(m-i+1)>0.
\eeqs Let $y=\frac {i+5}m\in (0.5, 1] $ and $ \ee=\frac 1m$. Then
\beqs \frac
{A}{8m}&=&6(1-y+2\ee)(y-\ee)(y-2\ee)(y-3\ee)(y-4\ee)-20(1-y+3\ee)(y-2\ee)(y-3\ee)(y-4\ee)\\&&
+30(1-y+4\ee)(y-3\ee)(y-4\ee) -30(1-y+5\ee)(y-4\ee)+15(1-y+6\ee)\\
&=&288\ee^5+(1584-744y)\ee^4+(720y^2-2340y+1920)\ee^3+(960+1270y^2-330y^3-1720y)\ee^2\\
&&+(210+72y^4-480y+510y^2-300y^3)\ee+15-45y-50y^3+60y^2+26y^4-6y^5.
\eeqs We can check that all these coefficients of $\ee$ are
nonnegative for $y\in (0.5, 1]$, so $A>0$ and the claim is proved.
\end{proof}

\begin{rem}The
sum of the last four terms $(A_0+A_1)(-\frac 2m)$  may be negative
when $\frac m2<i\leq m-2.$ In fact,  if $i=\frac {9}{10}m$ and $m$
is sufficiently large,  then  \beqs \frac {6}{i!}(A_0+A_1)(-\frac
2m)&=&-(m-i-2)(i+3)(i+2)(i+1)\frac
8{m^3}+(m-i-1)(i+2)(i+1)\frac {12}{m^2}\\
&&-(m-i)(i+1)\frac {12}m+6(m-i+1)\\&\sim&-0.0912m<0. \eeqs

\end{rem}

Next we consider the case $m-4\leq i\leq m-2$.

\begin{claim}\label{claim5}The lemma holds for $m-4\leq i\leq m-2$.
\end{claim}
\begin{proof}The proof is easy. If $i=m-4,$ then
\beqs \frac {6}{i!}P^{(m-4)}(-\frac 2m)&\geq&30-\frac
{48(m-3)}{m}+\frac {36}{m^2}(m-2)(m-3)-\frac
{16}{m^3}(m-1)(m-2)(m-3)\\&\geq &30-\frac {48(m-3)}{m}+\frac
{20}{m^2}(m-2)(m-3)\\
&=&\frac {2m^2+44m+120}{m^2}>0. \eeqs Similarly, we can prove that
the lemma holds for $i=m-3, m-2.$
\end{proof}

By Claim \ref{claim1}-\ref{claim5}, the lemma is proved.

\end{proof}

School of Mathematical Sciences, Peking University,  Beijing,
100871, P.R. China\\
Email: lihaozhao@gmail.com

\vskip3mm

\end{document}